\documentclass[11pt]{article}

\usepackage{amssymb}
\usepackage{latexsym,amsmath,amsthm,amsfonts}
\usepackage{graphicx, float, multicol, multirow,rotating}
\usepackage[dvips]{epsfig}
\usepackage{epstopdf} \usepackage{indentfirst}

\newtheorem{prop}{Proposition}
\newtheorem{lem}{Lemma}
\newtheorem{thm}{Theorem}

\newcommand{\be}{\begin{equation}}
\newcommand{\ee}{\end{equation}}

\newcommand{\bea}{\begin{eqnarray}}
\newcommand{\eea}{\end{eqnarray}}

\def\var{\mathop{\rm Var}}

\def\bW{{\bf W}}

\begin{document}

\bibliographystyle{plain}

\title{On the maximal size of Large-Average and ANOVA-fit Submatrices
in a Gaussian Random Matrix}

\author{Xing Sun and Andrew B.\ Nobel
\thanks{
Xing Sun is with Merck \& Co., Inc., One Merck Drive, Whitehouse Station,
New Jersey, 08889.  Email: xing\_sun$@$merk.com.
Andrew Nobel is with the Department of Statistics and Operations Research,
University of North Carolina, Chapel Hill, NC 27599-3260.
Email: nobel$@$email.unc.edu }}

\date{June 09, 2010}

\maketitle

\begin{abstract}

We investigate the maximal size of distinguished submatrices of a Gaussian
random matrix.  Of interest are submatrices whose entries have average
greater than or equal to a positive constant, and submatrices whose
entries are well-fit by a two-way ANOVA model.
We identify size thresholds and associated (asymptotic) probability bounds
for both large-average and ANOVA-fit submatrices.
Results are obtained when the matrix and submatrices of interest
are square, and in rectangular cases when the matrix submatrices
of interest have fixed aspect ratios.
In addition, we obtain a strong, interval concentration result
for the size of large average submatrices in the square case.
A simulation study shows good
agreement between the observed and predicted sizes of
large average submatrices in matrices of moderate size.
\end{abstract}

\vskip.2in

\noindent
Running title: Maximal submatrices of a Gaussian random matrix

\vskip.2in


\newpage

\section{Introduction}

Gaussian random matrices (GRMs) have been a fixture in the
application and theory of multivariate analysis for many
years.  Recent
work in the field of random matrix theory has provided a wealth of
information about the eigenvalues
and eigenvectors of Gaussian, and more general, random
matrices.  Motivated by problems of data mining and the
exploratory analysis of large data sets, this paper considers
a different problem, namely
the maximal size of distinguished submatrices in a GRM.
Of interest are submatrices
that are distinguished in one of two ways:
(i) the average of their entries is greater than or equal to a positive constant
or (ii) the optimal two-way ANOVA fit of their entries
has average squared residual less than a positive constant.

Using arguments from combinatorial probability, we
identify size thresholds and associated probability bounds
for large average and ANOVA-fit submatrices.  Results are obtained when
the matrix and the submatrices of interest are square, and when the matrix
and the submatrices of interest have fixed aspect ratios.
In each case, the maximal size of a distinguished submatrix grows
logarithmically with the dimension of the matrix, and depends in a
polynomial-type fashion on the inverse of the constant that constitutes
the threshold of distinguishability.  In the rectangular case,
the aspect ratio of the submatrix plays a more critical
role than the aspect ratio of the matrix itself.
In addition, we obtain upper and lower bounds
for the size of large average submatrices in the square case.
In particular, for $n \times n$ GRMs, the size
of the largest square submatrix with average greater than
$\tau > 0$ is eventually almost surely within in an interval
of fixed width that contains the critical value
$4 \tau^{-2} (\ln n -  \ln (4 \tau^{-2} \ln n))$.

We assess our bounds for large average submatrices
via a simulation study in which the size thresholds for
large average submatrices are compared to the observed
size of such submatrices in a Gaussian random matrix.
For matrices with moderate size and aspect ratio,
there is good agreement between the observed and predicted
sizes.

Results of the sort established here fall outside the purview of
random matrix theory and its techniques.   Nevertheless,
random matrix theory does provide some
insight into the logarithmic scale of large average submatrices.
This is discussed briefly in Section \ref{SVD} below.

\subsection{Exploratory Data Analysis}
\label{EDA}

The results of this paper are motivated in part by the increasing
application of exploratory tools such as biclustering to the
analysis of large data sets.  To be
specific, consider an $m \times n$ data matrix $X$
that is generated by measuring
the values of $m$ real-valued variables on each of $n$ subjects or samples.
The initial analysis of such data often
involves an exploratory search for interactions among samples
and variables.  In genomic studies of cancer,
sample-variable interactions can provide the basis for new insights
and hypotheses concerning disease subtypes and genetic pathways,
{\em c.f.}\
\cite{golub1999mcc, perou2000mph, garber2001dge, sorlie2001gep, sorlie2003rob, weigelt2005mpa}.

Formally, sample-variable interactions correspond to distinguished
submatrices of $X$.  The task of identifying such
submatrices is generally referred to as {\em biclustering}, two-way clustering or
subspace clustering in the computer science and bioinformatics
literature.  There is presently a substantial body of work on biclustering
methods, based on a variety of submatrix criteria; overviews can be found in
\cite{madeira2004bab, jiang2004cag, parsons2004sch}
and the references therein.  In particular, the biclustering methods by
Tanay \emph{et al.} \cite{tany1} and by Shabalin \emph{et al.} \cite{ShabalinandNobel} search
for submatrices whose entries have a large average value, while those of
Cheng and Church \cite{cheng00} and Lazzeroni and Owen \cite{owen02} search for submatrices
whose entries are well fit by a two-way ANOVA model.  The effectiveness
of these procedures in the analysis of real data is considered in
\cite{ShabalinandNobel}.


An exact or heuristic search among the (exponentially large)
family of submatrices of a
data matrix for those that are distinguished
by their average or ANOVA fit leads naturally
to a number of statistical questions related to multiple
testing.  For example, how large does a distinguished submatrix
have to be in order for it to be considered statistically significant, and
therefore potentially worthy of scientific interest?  What is the statistical
significance of a given distinguished submatrix?
Quantitative answers require an appropriate null
model for the observed data matrix, and in many cases, a GRM
model is a natural starting point for analysis.
When a GRM null is appropriate, the results of this paper provide partial
answers to the questions above.

We note that answers to statistical questions like those above
can have algorithmic implications.  For example,
knowing the minimal size of a significant submatrix can
provide a useful filtering criterion for exhaustive or heuristic search procedures,
or can drive the search procedure in a direct way.
The biclustering method in \cite{ShabalinandNobel} is based on a simple, Bonferroni
corrected measure of statistical significance that arises in the initial analyses
below.

\subsection{Bipartite Graphs}
\label{BG}

Our results on large average submatrices can also be
expressed in graph-theoretic terms, as every $m \times n$
matrix $X$ is associated in a natural way with
a bipartite graph $G=(V,E)$.  In particular, the vertex set
$V$ of $G$ is the disjoint union of two sets $V_1$ and $V_2$,
with $|V_1| = m$ and $|V_2| = n$, corresponding to the rows and
columns of $X$, respectively.
For each row $i \in V_1$ and column $j \in V_2$
there is an edge $(i,j) \in E$ with weight $x_{i,j}$.
There are no edges between vertices in $V_1$ or
between vertices in $V_2$.
With this association, large average submatrices of $X$
are in 1:1 correspondence with subgraphs of $G$ having
large average edge-weight.  The complexity of
finding the largest subgraph of $G$
whose average edge weight is
greater than a threshold appears to be unknown.  However,
it is shown in \cite{dhi1} that a slight variation of this problem,
namely finding the maximum edge weight subgraph in a general
bipartite matrix, is NP-complete.  A randomized, polynomial time
algorithm that finds a subgraph whose edge weight
is within a constant factor of the
optimum is described in \cite{alon04}, but this algorithm cannot
readily be adapted to the problem considered here.

\subsection{Connections with Random Matrix Theory}
\label{SVD}

The theory of random matrices provides some insight
into the relationship between
large average submatrices and the singular value decomposition.
In practice, the GRM assumption
made here acts as a null hypothesis.  If
an observed matrix contains a large average submatrix
whose size exceeds the thresholds given below, one may
reject the GRM hypothesis, and subject the identified
submatrix to further analysis.
This suggest an alternative hypothesis, under which a
fixed constant is added to every element of a select submatrix
of the null matrix, effectively embedding a large average submatrix
within a background of Gaussian noise.  It is then natural to ask if
the embedded submatrix affects the top singular
value or singular vectors of the resulting matrix.  We argue below
that the answer is a qualified no.

Let $W$ be an $m \times n$ Gaussian random matrix, representing
the null distribution.  Define a rank-one matrix
$S = 2 \tau a b^t$, where $\tau > 0$ is a fixed constant, and
$a \in \{0,1\}^m$, $b \in \{0,1\}^n$
are indicator vectors having $k$ and $l$ non-zero
components, respectively.  The outer produce $a \, b^t$
defines a submatrix $C$ whose rows
and columns are indexed by the indicator vectors $a$ and $b$, respectively.
The matrix $Y = W + S$ is distributed according to
an alternative hypothesis under which the fixed constant $\tau$ has
been added to every entry of the submatrix $C$.

Suppose that the dimensions $m,n,k$ and $l$ grow (with $n$, say)
in such a way that the matrix aspect ratio
$m / n \to \alpha$ with $\alpha \in [1,\infty)$,
and the submatrix aspect ratio $k / l$
remains bounded away from zero and infinity.
It is easy to see that the average of the $k \times l$ submatrix
$C$ in $Y$ has distribution
${\cal N}(2 \tau, (kl)^{-1})$, which is greater than $\tau$ with
overwhelming probability when $k$ and $l$ are large.
It follows from Proposition \ref{sig-bnd-avg} that
the probability of finding a $k \times l$ submatrix with average
greater than $\tau$ in the matrix $W$ is vanishingly small if $k$ and $l$
grow faster than $\log n$.  Thus, we might expect to see evidence of
$C$ in the first singular value, or the associated singular vectors, of $Y$.

Given an $m \times n$ matrix $U$,
let $s_1(U) \geq \cdots \geq s_m(U)$ denote its
ordered singular values, and let
$|| U ||_F = \sum_{i,j} u_{i,j}^2$ denote its Frobenius norm.
The difference between the largest singular value of
$W$ and $Y$ can be bounded as follows:
\begin{eqnarray}
\label{lidskii}
( s_1(Y) - s_1(W) )^2
& \leq &
\sum_{j=1}^n ( s_j(Y) - s_j(W) )^2  \nonumber \\
& \leq &
\sum_{j=1}^n (s_j(Y-W) )^2 \nonumber \\
& = &
|| Y - W ||_F^2 \ = \ || Z ||_F^2 \ = \ \tau^2 \, k \, l.
\end{eqnarray}
The second line above follows an inequality of
of Lidskii ({\em c.f.} Exercise 3.5.18 of \cite{matrixhj}), and the third makes
use of the fact that the Frobenius norm of a matrix is the
sum of the squares of its singular values.
By a basic result of Geman \cite{geman},
\be
\label{geman}
\frac{s_1(W)}{n^{1/2}} \to \left( 1 + \alpha^{1/2} \right)
\ee
with probability one as $n$ tends to infinity.
If $k = o(m^{1/2})$ and $l = o(n^{1/2})$, inequality
(\ref{lidskii}) implies that
$n^{-1/2} | s_1(Y) - s_1(W) | \to 0$ with probability one,
and therefore (\ref{geman}) holds with $Y$ in place of $W$.
In other words, the asymptotic behavior of $n^{-1/2} s_1(W)$ is
unchanged under the alternative $Y = W + Z$ if the
dimensions of the embedded submatrix $C$ grow more slowly
than $n^{1/2}$.  (Recall that $m$ is asymptotically proportional
to $n$.)

For fixed $\tau$ and $k,l$ such that $\log n < < k, l < < n^{1/2}$, the
embedded submatrix $C$ in $Y$ is highly significant,
but has no effect on the scaled limit of $s_1(Y)$.
Under the same conditions, $C$ is also not recoverable from the
top singular vectors of $Y$.
To be precise, let $u_1$ and $v_1$ be the left and right singular
vectors of $Y$ corresponding to the maximum singular
value $s_1(Y)$.  Using results of Paul \cite{paul} on the
singular vectors of spiked population models, it can be shown
that $a^t u_1$ and $b^t v_1$ tend to zero in probability
as $n$ tends to infinity.
Thus the row and column index vectors of
$C$ are asymptotically orthogonal to the first left and right
singular vectors of $Y$.

\subsection{Overview}

The next section contains probability bounds and a finite interval concentration
result for the size of large average submatrices in the square case.
Size thresholds and probability bounds for ANOVA submatrices in the square
case are presented in Section \ref{ANOVA}.  Thresholds and bounds in the rectangular
case are given in Section \ref{NSS}.  Section \ref{SIM} contains a
simulation study for large average submatrices.
Sections \ref{PLP} -- \ref{Basicpf}
contain the proofs of the main results.

\section{Thresholds and Bounds for Large Average Submatrices}
\label{SBAC}

Let $W = \{ w_{i,j} : i,j \geq 1 \}$ be an infinite array of independent ${\cal N}(0,1)$
random variables, and for $n \geq 1$, let $W_n = \{ w_{i,j} : 1 \leq  i,j \leq n \}$ be
the $n \times n$ Gaussian random matrix equal to upper left hand corner of $W$.
(The almost-sure asymptotics of Theorem \ref{Basic} requires
consideration of matrices $W_n$ that are derived from a fixed, infinite array.)
A submatrix of $W_n$ is a collection
$U = \{ w_{i,j} : i \in A, j \in  B \}$ where $A,B \subseteq \{1,\ldots,n\}$.
The Cartesian product $C = A \times B$ will be called the index
set of $U$, and we will write $U = W_n[C]$.  The dimension of $C$ is
$|A| \times |B|$, where $|A|, |B|$ denote the cardinality of $A$ and $B$,
respectively.  Note that the rows $A$ need not be contiguous,
and that the same is true of the columns $B$.
When no ambiguity will
arise, the index set $C$ will also be referred to as a submatrix of $W_n$.

\vskip.1in

\noindent
{\bf Definition:} For any submatrix $U$ of $W_n$ with index set $C = A \times B$, let
\[
F(U)
\  =  \
\frac{1}{|C|} \sum_{(i,j) \in C} w_{i,j}
\ = \
\frac{1}{|A| |B|} \sum_{i \in A, j \in B} w_{i,j}
\]
be the average of the entries of $U$.
Note that $F(U) \sim {\cal N}(0, |C|^{-1})$.

\vskip.1in

We are interested in the maximal size of square
submatrices whose averages exceed a fixed threshold.  This motivates
the following definition.

\vskip.1in

\noindent
{\bf Definition:} Fix $\tau>0$ and $n \geq 1$.  Let $K_{\tau}(W_n)$ be the largest $k \geq 0$
such that $W_n$ contains a $k \times k$ submatrix $U$ with $F(U) \geq \tau$.

\vskip.1in

As the rows and columns of a submatrix need not be contiguous, the statistic
$K_{\tau}(W_n)$ is invariant under row and column permutations of $W_n$.
We may regard the Gaussian distribution of $W_n$ as a null hypothesis for testing
an observed $n \times n$ data matrix, and
$K_{\tau}(\cdot)$ as a test statistic with which we can detect departures from
the null.
Our immediate goal is to obtain bounds on
the probability that $K_{\tau}(W_n)$ exceeds a given threshold, and to
identify a threshold for $K_{\tau}(W_n)$ that governs its asymptotic behavior.
To this end, we begin the analysis of $K_{\tau}(W_n)$ using standard first moment
type arguments, which are detailed below.

Let $\Gamma_k(n,\tau)$ be the number of $k \times k$ submatrices in
$W_n$ having average greater than or equal to $\tau$.
We begin by identifying the value of
$k$ for which $E \, \Gamma_k(n,\tau)$ is approximately equal to one.
If ${\cal S}_k$
denotes the set of all $k \times k$ submatrices of $W_n$ then
\be
\label{ukdef}
\Gamma_k(n,\tau)
\ = \
\sum_{U \in {\cal S}_k} \, I\{ F(W_n[U]) \geq \tau \} ,
\ee
and consequently
\be
\label{eukbnd}
E \, \Gamma_k(n,\tau)
\ = \ |{\cal S}_k| \cdot P(F(W_n[U]) \geq \tau)
\ = \ {n \choose k}^2 (1 - \Phi(\tau k) )
\ \leq \ {n \choose k}^2 e^{-\frac{\tau^2 k^2}{2}} ,
\ee
where in the last step we have used a standard bound on $1 - \Phi(\cdot)$.
For $s \in (0,n)$ define
\be
\label{phidef}
\phi_{n,\tau}(s)
\ = \
(2\pi)^{-\frac{1}{2}} \, n^{n+\frac{1}{2}} \, s^{-s-\frac{1}{2}} \,
(n-s)^{-(n-s)-\frac{1}{2}} \, e^{-\frac{\tau^2 s^2}{4}} .
\ \ \ \
\ee
Using the Stirling approximation of ${n \choose k}$, it is easy to see
that $\phi_{n,\tau}(k)$ is an approximation of the square root of the
final expression in (\ref{eukbnd}).   In particular, the
rightmost expression in (\ref{eukbnd}) is less than $2 \phi_{n,\tau}(k)^2$.
With this in mind, let $s(n,\tau)$ be any positive, real
root of the equation
\be
\label{sntaudef}
\phi_{n,\tau}(s) = 1 .
\ee
The next result shows that $s(n,\tau)$ exists and is unique, and it provides an
explicit expression for its value when $\tau$ is fixed and $n$ is large.

\begin{lem}
\label{sntau}
Let $\tau > 0$ be fixed.
When $n$ is sufficiently large, equation (\ref{sntaudef}) has a
unique root $s(n,\tau)$, and
\be
\label{k1n2}
s(n, \tau)
\ = \
\frac{4}{\tau^2} \, \ln n \, - \,
\frac{4}{\tau^2} \, \ln \left( \frac{4}{\tau^2} \ln n \right) \,
+ \, \frac{4}{\tau^2} \, + \, o(1)
\ee
where $o(1) \to 0$ as $n \to \infty$.
\end{lem}

\vskip.1in

We show below that the asymptotic behavior of the random
variables $K_{\tau}(W_n)$ is governed by the root $s(n,\tau)$
of equation (\ref{sntaudef}).
To begin, note that
for values of $k$ greater than $s(n,\tau)$, the expected number of
$k \times k$ submatrices $U$ of $W_n$ with $F(U) \geq \tau$ is
less than one.
The next proposition shows that the probability of seeing
such large submatrices is small.

\begin{prop}
\label{sig-bnd-avg}
Let $\tau>0$ be fixed.  For every $\epsilon > 0$,
when $n$ is sufficiently large,
\be
\label{sig-bnd1}
P\left( K_\tau(W_n) \geq  s(n,\tau) + r \right)
\ \leq \
\frac{4}{\tau^{2}} \, n^{- 2 \, r } \, \left( \frac{\ln n}{\tau^2} \right)^{2r + \epsilon}
\ee
for every $r = 1, \ldots, n$.
\end{prop}

\vskip.1in

The proofs of Lemma \ref{sntau} and Proposition \ref{sig-bnd-avg}
are given in Section \ref{PLP}.
The arguments are similar to those in \cite{sn08}, with adaptations to the
present setting.  A result similar to Proposition \ref{sig-bnd-avg}
can also be obtained from
the comparison principle for Gaussian sequences ({\em cf.}\ \cite{slep62}).
To be specific, fix $k \geq 1$ and note that
the family of random variables $\{ F(U) : U \in {\cal S}_k \}$
is a Gaussian
random field with $m = {n \choose k}^2$ elements that are pairwise
positively correlated, and have a common ${\cal N}(0,k \tau)$
distribution.
Then, by the comparison principle,
\[
P\left(K_{\tau}(W_n) \geq k \right)
\ =\
P \left( \max_{U \in {\cal S}_k} F(U) \geq \tau \right)
\ \le \
P\left(\max\{Z_1,...,Z_m\} \geq \tau \right),
\]
where $Z_1,\ldots,Z_m$ are independent ${\cal N}(0,k \tau)$ random
variables.  Using Poisson approximation based bounds such as those
in Section 4.4 of \cite{gold90}, one may obtain a
probability upper
bound similar to that in (\ref{sig-bnd1}).

\vskip.1in

It follows from Proposition \ref{sig-bnd-avg} and the Borel
Cantelli Lemma that, with probability one, $K_\tau(W_n)$ is eventually
less than or equal to $\lceil s(n,\tau) \rceil + 1 \leq s(n,\tau) + 2$.
Our principal result, stated in Theorem \ref{Basic} below, makes use of
a second moment argument in order to obtain
a corresponding lower bound.
The proof is given in Section \ref{Basicpf}.

\begin{thm}
\label{Basic}
Let $W_n$, $n \geq 1$, be Gaussian random matrices derived from an
infinite array $W$, and let $\tau > 0$ be fixed.  With probability one,
when $n$ is sufficiently large,
\be
\label{thm1}
s(n,\tau) \, - \, \frac{4}{\tau^2} \, - \, \frac{12 \ln 2}{\tau^2}
\, - \, 4
\ \leq \
K_{\tau}(\bW_n)
\ \leq \
s(n,\tau) + 2.
\ee
\end{thm}

The difference between the upper and lower
bounds in Theorem \ref{Basic} is a constant that depends on $\tau$, but
is {\em independent} of the matrix dimension $n$.  In particular the values
of the random variable $K_{\tau}(W_n)$ are eventually
concentrated on an interval that contains
$s(n,\tau)$ and whose width is independent of $n$.

The lower bound in Theorem \ref{Basic} can be further improved.
An examination of the argument in Lemma \ref{gauss2} in the Appendix
shows the inequality of the theorem still holds if the quantity
$12 \ln 2$ is replaced with any constant greater than $8 \ln 2$.

Extending earlier work of Dawande \emph{et al.} \cite{dhi1} and Koyuturk \emph{et al.} \cite{Gar2},
Sun and Nobel \cite{sn06,sn08} obtained a similar, two-point concentration
result for the size of
largest square submatrix of ones in an i.i.d.\ Bernoulli random matrix.
Bollob\'as and Erd\H{o}s \cite{bol2}, and Matula \cite{Mat1}, established analogous
results for the clique number of a regular random graph.
(See \cite{sn08} for additional references to work in the binary case.)
The proof of Theorem \ref{Basic} relies on a second moment argument,
but differs from the proofs of these earlier results
due to the continuous setting.
In particular, the proof makes use of the fact that,
under the Gaussian assumption made here,
for any $k \times k$ submatrix $U$ of $W$,
there exist simple upper bound and lower bounds on $P( F(U) \geq \tau )$,
and that the ratio of these bounds is of order $\tau k$.

\section{Thresholds and Bounds for ANOVA Submatrices}
\label{ANOVA}

In this section we derive bounds like those in Proposition \ref{sig-bnd-avg} for the size
of submatrices whose entries are well-fit by a two-way ANOVA model.
Roughly speaking, the ANOVA criterion
identifies submatrices whose rows (and columns) are shifts of one another.

\vskip.1in

\noindent
{\bf Definition:}
For a submatrix $U$ of $W_n$ with index set $A \times B$, define
\[
G(U)
\ = \
\min \left\{
\frac{1}{(|A|-1)(|B|-1)}
\sum_{i \in A, j \in B} (w_{ij} - a_i - b_j - c)^2
\right\} ,
\]
where the minimum is taken over all real constants $\{ a_i : i \in A \}$,
$\{ b_j : j \in B \}$ and $c$.

\vskip.1in

Under the ANOVA criterion, a submatrix $U$ will warrant interest
if $g(U)$ is less than a pre-defined threshold.
Note that by standard arguments,
\[
G(U)
\ = \
\frac{1}{(|A|-1)(|B|-1)}
\sum_{i \in A, j \in B} (w_{ij} - \overline w_{i.} - \overline w_{.j} + \overline w_{..})^2,
\]
where $\overline w_{i.}$, $\overline w_{.j}$, and $\overline w_{..}$ denote
the row, column, and the full submatrix averages, respectively.

\vskip.1in

\noindent
{\bf Definition:} Given $0< \tau < 1$, let $L_{\tau}(W_n)$ be
the largest value of $k$ such that $W_n$ contains a $k \times k$
submatrix $U$ with $G(U) \leq \tau$.

\vskip.1in

Arguments similar to those in the proof of Proposition \ref{sig-bnd-avg}, in
conjunction with
a probability upper bound on the
left tail of a $\chi^2$ distribution, establish the following
bound on $L_{\tau}(W_n)$.  The
proof is given in Section \ref{anovapf}.

\begin{prop}
\label{sig-bnd-anova}
Let $\tau > 0$ be fixed.  For every $\epsilon > 0$,
when $n$ is sufficiently large,
\be
\label{anova-ineq}
P\left( L_{\tau} (W_n) \geq t(n,\tau) + r \right)
\ \leq \
\frac{4}{h(\tau)}
\left(\frac{\ln n}{h(\tau)}\right)^{2r+2+\epsilon} n^{- 2 \, r }
\ee
for every $r = 1, \ldots, n$,
where
\[
t(n,\tau) = \frac{4}{h(\tau)}\ln n-\frac{4}{h(\tau)}
\ln \left(\frac{4}{h(\tau)}\ln n\right) + \frac{4}{h(\tau)} + 2
\]
and
\be
\label{hdef}
h(\tau)=1-\tau-\log(2-\tau).
\ee

\end{prop}

\vskip.1in

Proposition \ref{sig-bnd-anova} and the Borel Cantelli Lemma
imply that $L_{\tau}(W_n) \leq t(n,\tau) +1$ eventually almost surely.
The arguments used to lower bound $K_\tau(W_n)$ in Theorem
\ref{Basic} do not extend readily to $L_\tau(W_n)$, and we are
not aware if a similar interval-concentration result holds in this case.

\section{Thresholds and Bounds for Rectangular Submatrices}
\label{NSS}

The probability bounds of Proposition \ref{sig-bnd-avg} and
\ref{sig-bnd-anova} can be extended to non-square submatrices of
non-square matrices by adapting the methods of proof detailed in
Sections \ref{PLP} and \ref{anovapf}, respectively.  We present
the resulting bounds below, without proof.  Similar results concerning
submatrices of 1s in binary matrices can be found in \cite{sn08}.

\vskip.1in

\noindent
{\bf Definition:}
Let $W(m,n)$ denote an $m \times n$ Gaussian random matrix, and let
$\alpha > 0$ and $\beta \geq 1$ be fixed aspect ratios for the sample
matrix and target submatrix respectively.

a. For $\tau > 0$ let $K_\tau(W : n, \alpha, \beta)$ be the largest integer
$k$ such that there
exists a $\lceil \beta k \rceil \times k$ submatrix $U$ in
$W(\lceil \alpha n \rceil, n)$ with $F(U) \geq \tau$.

b. For $0 < \tau < 1$ let $L_\tau(W : n, \alpha, \beta)$
be the largest integer k such that there
exists a $\lceil \beta k \rceil \times k$ submatrix $U$ in
$W(\lceil \alpha n \rceil, n)$ with $G(U) \leq \tau$.

\vskip.2in

\begin{prop}
\label{non-sq-avg} Fix $\tau>0$ and any $\epsilon >0$. When $n$ is sufficiently
large,
\[
P\left( K_\tau(W : n, \alpha, \beta) \ \geq \ s(n, \tau, \alpha, \beta) + r \right)
\ \leq \
n^{- (\beta+1) \, r}
\left( \frac{\ln n}{\tau^2} \right)^{(\beta+1+\epsilon) r}
\]
for each $1 \le r \le n$, where
\[
s(n, \tau, \alpha, \beta)
\ = \
\frac{2(1 + \beta^{-1})}{\tau^2} \ln n -
\frac{2(1 + \beta^{-1})}{\tau^2} \ln \left[ \frac{2(1 + \beta^{-1})}{\tau^2} \ln n \right]
+ \frac{2}{\tau^2} \ln \alpha+C_1(\beta,\tau),
\]
for some constant $C_1(\beta,\tau)>0$.
\end{prop}

\vskip.2in

\begin{prop}
\label{non-sq-anova} Fix $0<\tau<1$ and any $\epsilon >0$. When $n$ is
sufficiently large,
\[
P( L_\tau(W : n, \alpha, \beta) \ \geq \ t(n, \tau, \alpha, \beta) + r ) \
\ \leq \
n^{- (\beta+1) \, r} \left(\frac{\ln n}{h(\tau)}\right)^{(\beta+1+\epsilon)r}
\]
for each $1 \le r \le n$, where
\[
t(n, \tau, \alpha, \beta)
\ = \
\frac{2(1 + \beta^{-1})}{h(\tau)} \ln n -
\frac{2(1 + \beta^{-1})}{h(\tau)} \ln \left[ \frac{2(1 + \beta^{-1})}{h(\tau)} \ln n \right] +
h(\tau)^{-1}\ln \alpha +C_2(\beta,\tau),
\]
for some constant $C_2(\beta,\tau)>0$, where $h(\tau)$ is defined as in (\ref{hdef}).
\end{prop}

\noindent
{\bf Remark:} The bounds in Propositions \ref{non-sq-avg}
and \ref{non-sq-anova} have a similar form.  In each case, the
bound is of the form $n^{- (\beta+1) \, r}$ times a polynomial
in $\ln n$, and the leading term in $s(\cdot)$
and $t(\cdot)$ are of the form $(1 + \beta^{-1}) \ln n$ times
a function of the threshold $\tau$.  We note the critical role
played by the aspect ratio $\beta$ of the target submatrix.
By contrast, the aspect ratio $\alpha$ of the sample matrix plays
a secondary role, its logarithm appearing only in the constant term of
$s(\cdot)$ and $t(\cdot)$.

\section{Simulation Study for Large Average Submatrices}
\label{SIM}

The size thresholds and probability bounds presented in
Sections \ref{SBAC} - \ref{NSS} are asymptotic, and it is
reasonable to ask if they apply to matrices of moderate
size.  To this end, we carried out a simulation study
in which we compared the
size of large average submatrices in simulated Gaussian
data matrices with the bounds predicted by the theory.
An exhaustive search for large average submatrices is
not computationally feasible.
Our study was based on a simple search algorithm
for large average submatrices that is used in the biclustering
procedure of Shabalin
{\it et al.} \cite{ShabalinandNobel}.  Analogous application of
existing ANOVA based biclustering procedures does not appear
to be straightforward, so the simulation study was restricted to
the large average criteria.

The search algorithm from \cite{ShabalinandNobel} operates
as follows.  Given an $m \times n$ data matrix $W$ and integers
$1 \leq k \leq m$ and $1 \leq l \leq n$, a random subset of $l$
columns of $W$ is selected.  The sum of each row over the
selected set of $l$ columns is computed, and the rows
corresponding to the $k$ largest sums are selected.
Then the sum of each column over the
selected set of $k$ rows is computed, and the columns
corresponding to the $l$ largest sums are selected.
This alternating update of row and column sets is repeated
until a fixed point is reached, and the average of the resulting
$k \times l$ matrix is recorded.  The basic search procedure is
repeated $N$ times, and the output of the search algorithm is the
largest of the $N$ observed submatrix averages.
The search algorithm is not guaranteed to find the $k \times l$
submatrix of $W$ with maximum average.  However, the algorithm
provides a lower bound on the maximum average value of
$k \times l$ submatrices
We conducted two experiments, one for square matrices
and one for rectangular matrices.

\vskip.1in

\begin{figure}
\caption{Results of 200 x 200 simulations}
\vskip.1in
\begin{center}
\resizebox{12cm}{6cm}{\includegraphics[angle=90]{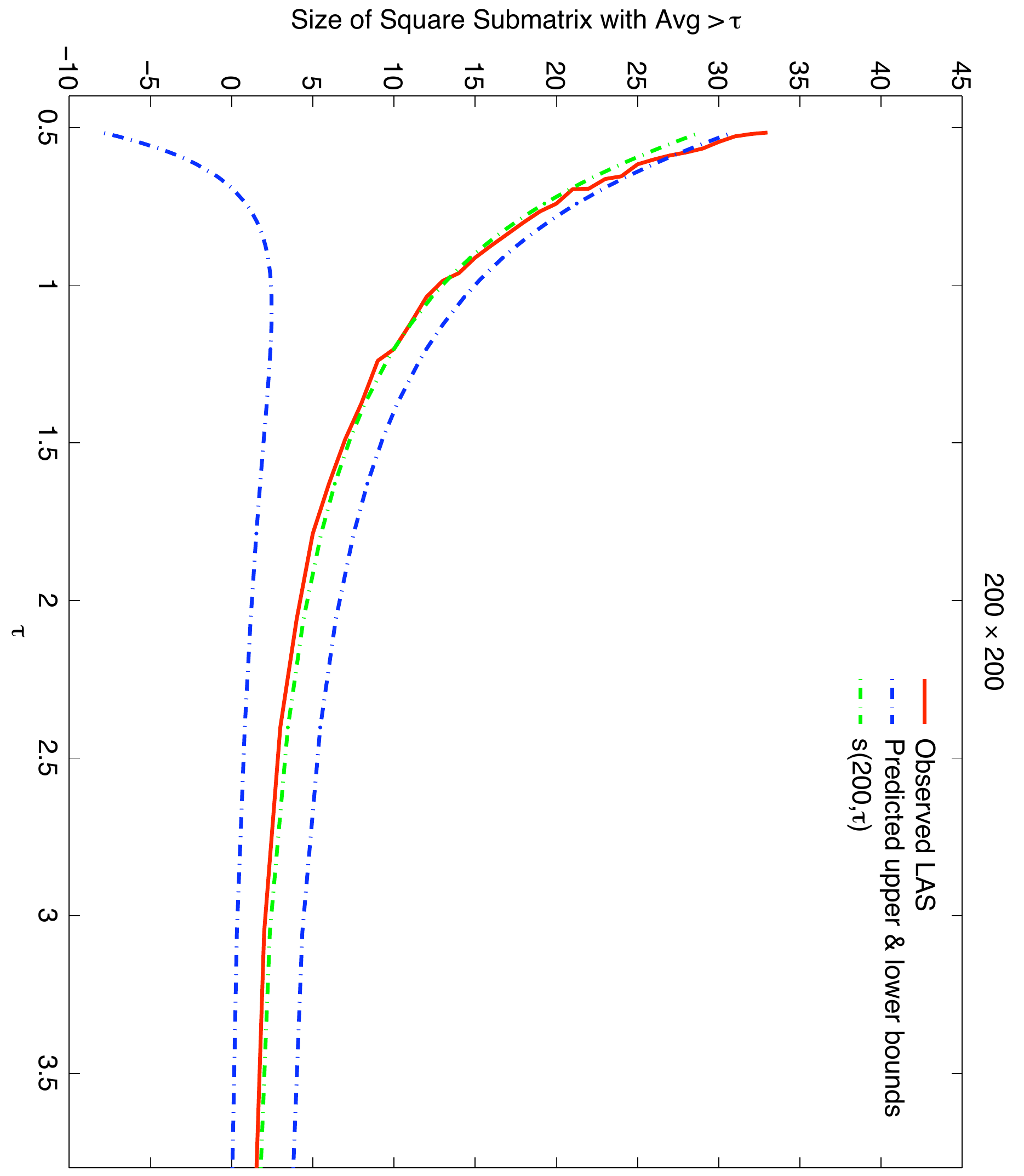}}
\end{center}
\label{sq200}
\end{figure}

{\bf Square matrices.}  We considered matrices of size
$n = 200$ and $n = 500$.
Results from the case $n= 200$ are summarized in Figure \ref{sq200}.
For a fixed $k \geq 1$, we generated a $200 \times 200$
Gaussian random matrix $W$, and then used the
search algorithm described above to find a lower
bound, $\tau_k$, on the maximum average of the
$k \times k$ submatrices of $W$ using $N = 10000$ iterations of the
search procedure.
Different random matrices $W$ were
generated for different values of $k$.
The upper and lower bounds of Theorem \ref{Basic} begin to diverge
when $\tau \leq 1/2$, so we restricted attention to values of $k$
for which $\tau_k > 1/2$.  In this case
$k$ ranged from 1 to 55.
A linear interpolation of the pairs $(\tau_k,k)$
appears as the red curve in Figure \ref{sq200}.  We have also plotted
the threshold function $s(n,\tau)$ derived in Lemma \ref{sntau}, omitting the
$o(1)$ term, as well as the upper and lower bounds from Theorem \ref{Basic}.
As can be seen from the figure, there is good agreement
between the observed and predicted sizes of large average submatrices.
In particular, for the range $\tau \geq 1/2$ the observed sizes of large
average submatrices fall within the upper and lower bounds of
the theorem.

\begin{figure}
\caption{Results of 500 x 500 simulations}
\vskip.1in
\begin{center}
\resizebox{12cm}{6cm}{\includegraphics[angle=90]{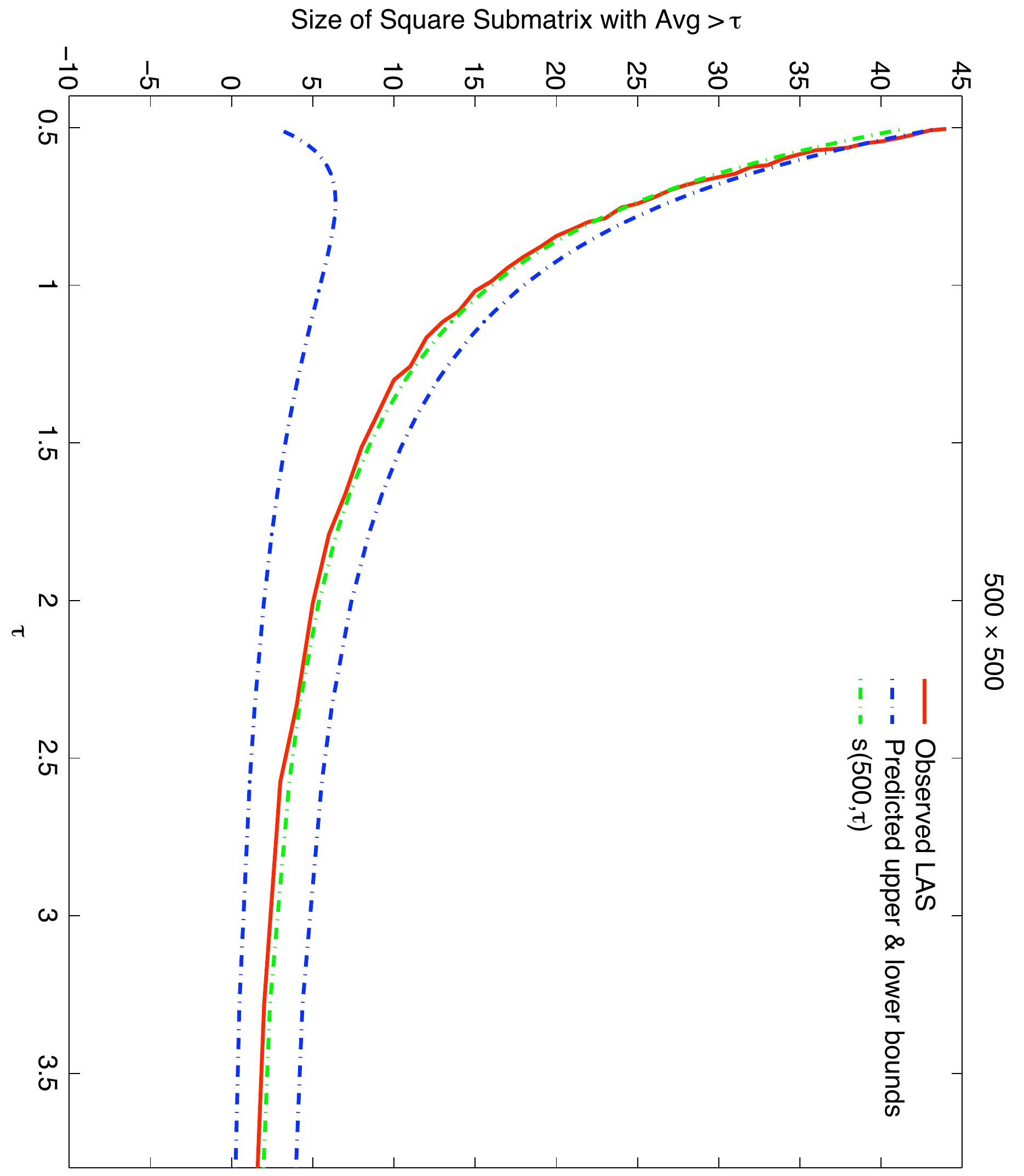}}
\end{center}
\label{sq500}
\end{figure}

Simulations for matrix size $n = 500$ were carried out in a similar
fashion.  The results, based on $N = 10000$ iterations of the search
procedure
for each value of $k$, are summarized in Figure \ref{sq500}.
Restricting attention to $\tau_k > 1/2$ leads to matrix
sizes $k$ between $1$ and $55$
As in the case $n = 200$
there is good agreement between the observed and predicted
sizes of large average submatrices, and the observed sizes of large
average submatrices fall within the upper and lower bounds of
Theorem \ref{Basic}.

\vskip.2in

\begin{figure}
\label{rcc}
\caption{Results for rectangular simulations}
\vskip.1in
\begin{center}
\resizebox{12cm}{6cm}{\includegraphics[angle=90]{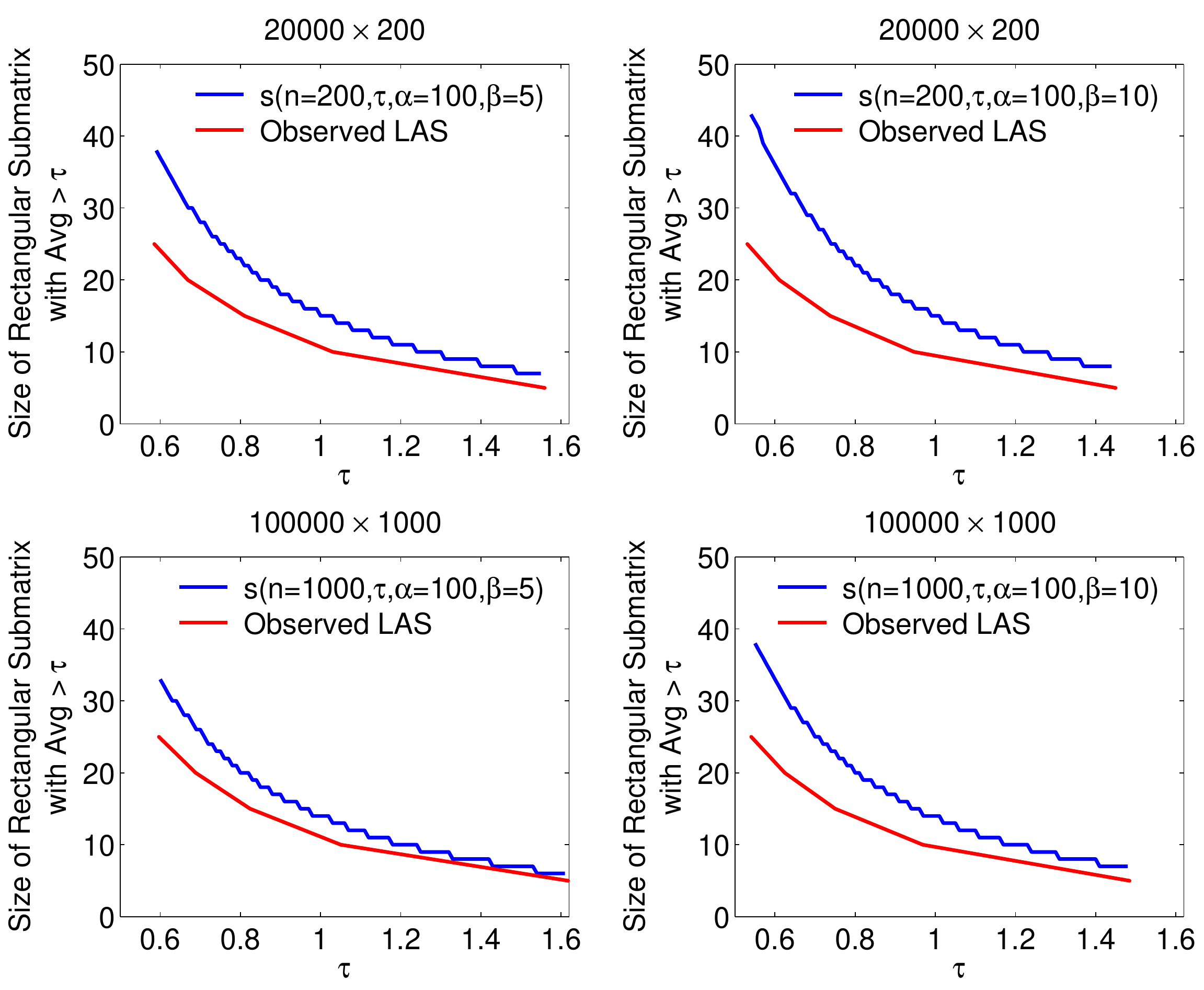}}
\end{center}
\end{figure}

{\bf Non-Square matrices.}
We also carried out two simulation studies for rectangular matrices
of sizes $20,000 \times 200$ and $100,000 \times 1000$
(matrix aspect ratio $\alpha=100$).
These sizes reflect those commonly seen in high-throughput
genomic data.  In each case, we looked for submatrices with
aspect ratio $\beta = 5$ and $\beta = 10$.
For each fixed $k \in \{5, 10, 15, 20, 25 \}$, we generated a
Gaussian random matrix of the appropriate size and then
used the search algorithm with $N = 10000$ iterations
to identify $\beta k \times k$ submatrices with large average.
The results are summarized in the (interpolated) red curves of Figure \ref{rcc}.
The theoretical upper bounds from Proposition \ref{non-sq-avg} are
plotted in blue for comparison.
In each case the observed maxima lie below the theoretical upper bound;
the gap decreases with decreasing $\beta$ and increasing $\tau$.

\section{Proof of Lemma \ref{sntau} and Proposition \ref{sig-bnd-avg}}
\label{PLP}

\noindent
{\bf Proof of Lemma \ref{sntau}:}
Let $\tau > 0$ be fixed, and note that
\be
\label{lnphi}
\ln \phi_{n,\tau}(s) \ =  \
(n+\frac{1}{2}) \ln n - (s+\frac{1}{2}) \ln s - (n-s+\frac{1}{2}) \ln (n-s)
 - \frac{\tau^{2} s^2}{4} - \frac{1}{2} \ln 2 \pi .
\ee
Differentiating $\ln \phi_{n,\tau}(s)$ with respect
to $s$ yields
\[
\frac{\partial \ln \phi_{n,\tau}(s)}{\partial s} \ = \
\frac{1}{2(n-s)} + \ln (n-s) - \frac{1}{2s} - \ln s - \frac{s \tau^2}{2} .
\]
The last expresssion is negative when
$2 \tau^{-2} \ln n < s < 4 \tau^{-2} \ln n$;
we now consider the value of $\ln \phi_{n,\tau}(s)$ for $s$ outside
this interval.
A straightforward
calculation shows that for $0 < s \leq 2 \tau^{-2} \ln n$,
\bea
\ln \phi_{n,\tau}(s)
& \geq &
s \, \left( \ln (n - 2\tau^{-2} \ln n) - \frac{s \tau^2}{4} - \ln \ln n -\ln 2\tau^{-2}\right)
- \frac{1}{2} \ln s  - \frac{1}{2} \ln 2 \pi , \nonumber
\eea
which is positive
when $n$ is sufficiently large.  In order to address the other
extreme, note that from (\ref{lnphi}) we have
\be
\label{lgineq}
\ln \phi_{n,\tau}(s)
\ \leq \
s \, \left( \ln (n - s) - \frac{s \tau^2}{4} - \ln s \right)
- \frac{1}{2} \ln s  + (n + 1/2) \ln \left( \frac{n}{n-s} \right) .
\ee
It is easy to check that the right hand side of the above inequality is negative
when $s>n-2$.
Considering separately the cases
$s+2 < n < (2 \ln 2)^{-1}{s\ln s}$
and $n \ge (2 \ln 2)^{-1}{s \ln s}$,
one may upper bound the final term above by
$(s \ln s)/2 +(\ln 2)/2$ and $2s+(\ln 2)/2$, respectively.
Thus, for $s < n-2$, we have
\[
\ln \phi_{n,\tau}(s)
\ \leq \
s \, \left( \ln (n - s) - \frac{s \tau^2}{4} - \ln s \right)
- \frac{1}{2} \ln s + 2s + \frac{ s \ln s }{2}
+ \frac{\ln 2}{2} ,
\]
and in particular, for $4 \tau^{-2} \ln n \leq s < n-2$,
\[
\ln \phi_{n,\tau}(s)
\ \leq \
s \, \left( 2 - \frac{\ln s}{2} \right)
- \frac{1}{2} \ln s  + \frac{\ln 2}{2}
\ < \ 0 \nonumber
\]
when $n$ (and therefore $s$) is sufficiently large.
Thus for large $n$ there exists a unique solution $s(n,\tau)$ of
the equation $\phi_{n,\tau}(s)=1$ with
$s(n,\tau) \in (2 \tau^{-2} \ln n, 4 \tau^{-2} \ln n)$.

Taking logarithms of both sides of the equation $\phi_{n,\tau}(s) = 1$ and rearranging terms
yields the expression
\be
\label{phieq1}
\frac{1}{2} \ln \frac{n}{n-s} + n \ln \frac{n}{n-s} -
(s+\frac{1}{2}) \ln s + s \ln (n-s) - \frac{\tau^2 s^2}{4}
\ =\
\frac{\ln 2 \pi}{2}.
\ee
The argument above shows that the (unique) solution of this
equation belongs to the interval $(2 \tau^{-2} \ln n, 4 \tau^{-2} \ln n)$, so
we consider the case in which $s$ and $n / s$ tend to infinity with $n$.
Dividing both sides of (\ref{phieq1}) by $s$ yields
\[
\ln (n-s) - \frac{s \tau^2}{4} - \ln s
\ = \
- 1 + O(\frac{\ln s}{s}),
\]
which, after adding and subtracting terms, can be rewritten in the
equivalent form
\be
\label{phieq2}
\ln n - \frac{s \tau^2}{4} - \ln \ln n
\ = \
\ln \left( \frac{s}{\ln n} \right) - \ln \left( \frac{n-s}{n} \right) - 1 + O(\frac{\ln s}{s}).
\ee
For each $n \geq 1$, define $R(n)$ via the equation
\[
s(n,\tau) \ = \ 4 \tau^{-2} \ln n- 4 \tau^{-2} \ln \ln n + R(n).
\]
Plugging the last expression into (\ref{phieq2}), we find that
$R(n) = \frac{4}{\tau^2}(1-\ln\frac{4}{\tau^2}) + o(1)$, and the result follows from
the uniqueness of $s(n,\tau)$.

\vskip.2in

\noindent
{\bf Proof of Proposition \ref{sig-bnd-avg}:}
Fix $\tau > 0$.
If $\lceil s(n,\tau) \rceil + r > n$ the bound
(\ref{sig-bnd-avg}) holds trivially; in the case of
equality, it follows from a standard Gaussian
tail bound when $n$ is sufficiently large.
Fix $n \geq 1$ for the moment and suppose that
$l = \lceil s(n,\tau) \rceil + r \leq n-1$.
By Markov's inequality and the definition of $\phi_{n,\tau}(\cdot)$,
\bea
P( M_\tau(W_n) \geq s(n,\tau) + r )
& = &
P( M_\tau(W_n) \geq l ) \nonumber \\
& = &
P(U_l(n,\tau) \geq 1) \nonumber \\
& \leq &
E U_l(n,\tau) \nonumber \\
& \leq &
\label{mark}
2 \phi_{n,\tau}^2(l)
\ \leq \
2 \phi_{n,\tau}^2( s(n,\tau) + r ) .
\eea
Let $\gamma = e^{- \tau^2 / 4}$ and, to reduce notation,
denote $s(n,\tau)$ by $s_n$.
Under the constraint on $r$, a straightforward calculation
shows that one can decompose the
final term above as follows:
\be
\label{phin2}
2 \phi_{n,\tau}^2(s_n + r)
\ = \
2 \, \phi_{n,\tau}^2(s_n) \, \gamma^{2r  s_n} \,
       [ \, A_n(r) \, B_n(r) \, C_n(r) \, D_n(r) \, ]^2
\ee
where
\bea
A_n(r)
& = &
\left( \frac{n - r - s_n}{n - s_n} \right)^{- n + r + s_n - \frac{1}{2}}
\ \ \ \ \ \ \
B_n(r)
\ = \
\left( \frac{r + s_n}{s_n} \right)^{- s_n - \frac{1}{2}} \nonumber \\
[.1in]
C_n(r)
& = &
\left( \frac{n - s_n}{r + s_n} \, \gamma^{s_n} \right)^r
\ \ \ \ \ D_n(r) \ = \ \gamma^{r^2}  \nonumber
\eea

\vskip.06in

It is enough to bound the right hand side of (\ref{phin2}) as $n$
increases and $r = r(n)$ is such that
$\lceil s(n,\tau) \rceil + r \leq n-1$.
By definition, $\phi_{n,\tau}(s_n) = 1$, and for each fixed $\epsilon>0$,
\[
\max_{r \geq 1} \,
\frac{ 2 \gamma^{2r s_n} }{ n^{-2r} (\frac{2\ln n}{\tau^2})^{2r+\epsilon} }
\to 0
\ \mbox{ as } \ n \to \infty.
\]
Thus it suffices to show that the product
$A_n(r) \, B_n(r) \, C_n(r) \, D_n(r)$ is uniformly bounded in $r$.
To begin, note that for any fixed
$0 < \delta < 4$,
\[
C_n(r)^{\frac{1}{r}}
\ = \
\frac{n - s_n}{r + s_n} \,
\gamma^{s_n}
\ \leq \
\frac{n}{s_n} \, \gamma^{s_n}
\ \leq \
\frac{4}{4 - \delta} \, e^{-1} \cdot o(1) .
\]
The last term will be less than one when
$\delta$ is sufficiently small.
The term $B_n(r) \leq 1$ for each $r \geq 1$, so it only
remains to show that $\max_{r \geq 1} A_n(r) \cdot D_n(r)$
is bounded as a function of $n$.  A straightforward calculation shows
that $\ln A_n(r) \leq r$, and consequently,
$\ln A_n(r) \cdot D_n(r) \le r - \frac{\tau^2r^2}{4}$, a quadratic
function of $r$ that is bounded from above.

\section{Proof of Proposition \ref{sig-bnd-anova}}
\label{anovapf}

For any $k \times k$ submatrix $U$ of the Gaussian random
matrix $W_n$, it follows from standard arguments that
$(k-1)^2 G(U)$ has a $\chi^2$
distribution with $(k-1)^2$ degrees of freedom.  In order to bound
the quantity $P(G(U) \leq \tau)$, which arises in the analysis
of $L_\tau(W_n)$, we
require an initial result relating the right and left tails of the
$\chi^2$ distribution.

\begin{lem}
\label{anova1}
Suppose that $X \sim \chi^2_\ell$ for some $\ell \ge 3$.
Then for $0 < t < \ell-2$ we have
\[
P(X \leq t) \ \leq \ P(X \geq 2\ell-4-t).
\]
\end{lem}

\noindent
{\bf Proof of Lemma \ref{anova1}:}
Let $f$ denote the density function of $X$ and let
$0 < t < \ell-2$.  Since
\[
P(X \le t)
\ = \
\int_0^{t} f(s) \, ds
\ \mbox{ and } \
P(X \geq 2\ell-4-t)
\ \geq \
\int_{2\ell-4-t}^{2\ell-4} f(s) \, ds,
\]
it suffices to show that
\be
\label{chif1}
\frac{f(s)}{f(2\ell-4-s)} \ \leq \ 1 \
\mbox{ for all }
0 < s < \ell-2 .
\ee
To this end, note that the ratio in (\ref{chif1}) can be rewritten as follows:
\bea
\frac{f(s)}{f(2\ell-4-s)}
& = &
\frac{s^{(\ell-2)/2} \, e^{-s/2}}
{(2\ell-4-s)^{(\ell-2)/2} \, e^{-(2\ell-4-s)/2}} \nonumber\\ [.08in]
& = &
\left[ \left(1-\frac{2\ell-4-2s}{2\ell-4-s} \right) e^{2(\ell-2-s) / (\ell-2)} \right]^{(\ell-2)/2} .
\nonumber \\ [.08in]
& = &
\left[ \left(1-\frac{1}{u} \right)  e^{\frac{2}{2u-1}} \right]^{(\ell-2)/2}
\ \mbox{ with } \ u = \frac{2\ell-4-s}{2\ell-4-2s} .
\label{monodiff}
\eea
As $s$ tends to $\ell-2$, $u$ tends to infinity, and therefore
\[
\lim_{s \to (\ell-2)}
\frac{f(s)}{f(2\ell-4-s)}
\ = \
\lim_{u \rightarrow \infty} \left( 1-\frac{1}{u} \right) e^{\frac{2}{2u-1}}
\ =\ 1.\]
Thus, it suffices to show that for $u \in(1,\infty)$, the final term
in (\ref{monodiff}) is an increasing function of $u$.
Differentiating with respect to $u$ we find that
\[
\frac{d}{du}
\left(1-\frac{1}{u} \right) e^{\frac{2}{2u-1}}
\ = \
\frac{(2u-1)^2-4(u-1)u}{u^2(2u-1)^2} \,
e^{\frac{2}{2u-1}}
\ > \ 0
\]
where the inequality follows from the fact that $u > 1$.
Inequality (\ref{chif1}) follows immediately.

\vskip.2in

\noindent
{\bf Proof of Proposition \ref{sig-bnd-anova}:}
To begin, note that if $X$ has a $\chi^2$ distribution with $\ell$ degrees
of freedom, then by a standard Chernoff bound,
\be
\label{Chernoff}
P( X \geq r)
\ \leq \
\min_{0 < s < \frac{1}{2}} (1-2s)^{-\frac{\ell}{2}} \, e^{-s r}
\ = \
\left[ \left( \frac{\ell}{r} \right)
e^{ \left(\frac{r}{\ell} - 1 \right) } \right]^{-\ell/2}
\ee

Let $\tau > 0$ be fixed.  Fix $n \geq 1$ for the moment
and let $r \geq 1$ be such that
$k = \lceil t(n,\tau) \rceil + r \leq n$, where
$t(n,\tau)$ is defined as in the statement of
Proposition \ref{sig-bnd-anova}.
Let $U$ be any $k \times k$ submatrix of $W_n$, and let
$\ell = (k-1)^2$.
As noted above, the random variable $\ell \, G(U)$ has a $\chi^2$
distribution with $\ell$ degrees of freedom, so by
Lemma \ref{anova1} and inequality (\ref{Chernoff}),
\bea
P(G(U) \leq \tau)
& = &
P( \ell \, G(U) \, \leq \, \ell\, \tau)
\ \leq \
P\left( \ell \, G(U) \, \geq \, (2 - \tau) \ell - 4 \right) \nonumber \\ [.1in]
& \leq &
\exp\left\{
- \frac{\ell}{2}
\left[  \frac{{(2 - \tau) \ell - 4}}{\ell} - 1 + \ln \frac{\ell}{(2 - \tau) \ell - 4} \right] \right\}
\nonumber \\[.1in]
& = &
\exp\left\{
- \frac{\ell}{2}
\left[  (1 - \tau) - \ln (2 - \tau) \right]
\right\} \,
\exp\left\{
\left[ 2 + \frac{\ell}{2} \ln \left(1 - \frac{4}{\ell (2 - \tau)} \right) \right]
\right\} .
\nonumber
\eea
One may readily show that the second term above is $O(1)$.
It then follows from a first moment argument that
\be
\label{ltau}
P( L_\tau (W_n) \geq k )
\ \leq \
{n \choose k}^2 P(G(U) \leq \tau)
\ \leq \
C \, {n \choose k}^2 q^{(k-1)^2}
\leq
C \, {n \choose k-1}^2 q^{(k-1)^2}\cdot n^2
\ee
where $C$ is a finite constant and
\[
q \ = \ \exp\left\{\frac{1}{2}[-(1-\tau) + \ln (2-\tau)]\right\} .
\]
Fix $\epsilon > 0$.  By following the proofs of
Lemma \ref{sntau} and Proposition \ref{sig-bnd-avg},
replacing $\tau^2$ with $h(\tau) = 1-\tau-\ln (2-\tau)$,
one can show that for every $r \geq 1$ such that
\[
k \ = \
\left\lceil \frac{4}{h(\tau)}\ln n-\frac{4}{h(\tau)}\ln \left(\frac{4}{h(\tau)}\ln n\right)
+ \frac{4}{h(\tau)} \right\rceil + 2 + r
\]
is at most $n$, we have
\[
{n \choose k-1}^2 q^{(k-1)^2}
\ \leq \
\frac{4}{h(\tau)}
\left(\frac{\ln n}{h(\tau)}\right)^{2r+2+\epsilon} n^{- 2 \, r - 2 } ,
\]
and the result then follows from (\ref{ltau}).

\section{Proof of Theorem \ref{Basic}}
\label{Basicpf}

In what follows we make use of standard bounds on the tails
of the Gaussian distribution, namely that
$(3s)^{-1} e^{-s^2 / 2} \leq 1 - \Phi(s) \leq s^{-1} e^{-s^2 / 2}$
for $s \geq 3$.  The proof of Theorem \ref{Basic} is based on
several preliminary results.  The first result bounds the ratio of
the variance of $\Gamma_k(\tau,n)$ and the square of its
expected value, a quantity that later arises from an application
of Chebyshev's inequality.

\begin{lem}
\label{gauss1}
Fix $\tau>0$.
There exist integers $n_0, k_0 \geq 1$ and a positive constant $C$ depending
on $\tau$ but independent of $k$ and $n$, such that
for any $n \ge n_0$ and any $k \ge k_0$,
\be
\frac{ \var \, \Gamma_k(\tau,n)}{(E \, \Gamma_k(\tau,n))^2}
\ \le \
C \, k^4 \,
\sum_{l=1}^k \sum_{r=1}^k
         \frac{{k \choose l}{n-k \choose k-l}}{{n \choose k}}
         \frac{{k \choose r}{n-k \choose k-r}}{{n \choose k}}
         \exp \left\{\frac{rl\tau^2}{2}\left(1+\frac{k^2-rl}{k^2+rl}\right)\right\}.
\label{bagt1}
\ee
\end{lem}

\noindent{\bf Proof:}
Let $\mathcal{S}_k$ denote the collection of all $k \times k$
submatrices of $W_n$. It is clear that
\be
E \, \Gamma_k(n,\tau)
=
\sum_{U \in \mathcal{S}_k} P(F(U) > \tau)
=
{n \choose k}^2 \left( 1-\Phi(k \tau) \right) .
\label{mean}
\ee
In a similar fashion, we have
\[
E \, \Gamma^2_k(n,\tau)
\ = \
\sum_{U_i, U_j \in \mathcal{S}_k} P( F(U_i)>\tau \mbox{ and } F(U_j) > \tau )
\]
Note that the joint probability in the last display depends only on the
overlap between the submatrices $U_i$ and $U_j$.
For $1 \leq r,l \leq k$ define
\[
G(r,l) \ = \ P( F(U) > \tau \mbox{ and } F(V) > \tau )
\]
where $U$ and $V$ are two fixed $k \times k$ submatrices of $W$ having
$r$ rows and $l$ columns in common.  Then
$E \, \Gamma_k(n,\tau) = {n \choose k}^2 G(0,0)^{1/2}$, and a straightforward
counting argument shows that
\[
E \, \Gamma^2_k(n,\tau)
\ = \
\sum_{r=0}^k
\sum_{l=0}^k
{n \choose k}^2
{k \choose r}{n-k \choose k-r}
{k \choose l}{n-k \choose k-l} \,
G(r,l) .
\]
In particular,
\bea
\frac{\mbox{Var} \, \Gamma_k(n,\tau)}{(E \, \Gamma_k(n,\tau))^2}
& = &
\sum_{r=0}^k \sum_{l=0}^k
\frac{ {k \choose l}{ n-k \choose k-l} }{ {n \choose k} }
\frac{ {k \choose r}{n-k \choose k-r} }{ {n \choose k} }
\left(\frac{ G(r,l) }{ G(0,0) }\right) - 1.
\nonumber \\
& = &
\sum_{r=1}^k \sum_{l=1}^k
\frac{ {k \choose l}{ n-k \choose k-l} }{ {n \choose k} }
\frac{ {k \choose r}{n-k \choose k-r} }{ {n \choose k} }
\left( \frac{ G(r,l) }{ G(0,0) } - 1 \right).
\nonumber
\eea
where we have used the fact that ${k \choose l}{n-k \choose k-l} / {n \choose k}$
is a probability mass function, and that $G(0,l) = G(r,0) = G(0,0)$.
When $k \tau \geq 3$ we have
$G(0,0) = (1- \phi(k \tau))^2 \geq (3 k \tau)^{-2} e^{- k^2 \tau^2}$,
and it therefore suffices to show that for $1 \leq r, l \leq k$,
\be
\label{grlineq}
G(r,l)
\ \leq \
C \, k^2 \,
\exp \left\{ -k^2 \tau^2  + \frac{ rl \tau^2 }{2} \left( 1 + \frac{k^2 - rl}{k^2 + rl} \right) \right\}
\ee
where $C > 0$ depends on $\tau$ but is independent of $k$ and $n$.
Inequality (\ref{grlineq}) is readily established when $r = l = k$, so we
turn our attention to bounding $G(r,l)$ when $1 \leq r l < k^2$.  In this
case
\[
G(r,l)
\ = \
\frac{\sqrt{rl}}{\sqrt{2\pi}} \int_{-\infty}^{\infty}  e^{-\frac{rlt^2}{2}} \,
P \left(F(U \cap V^c) \geq \frac{k^2 \tau - rlt}{\sqrt{k^2 - rl}} \right)^2 \, dt
\]
where $U,V$ are submatrices of $W_n$ having $r$ rows and
$l$ columns in common.  Let $\overline{\Phi}(x) = 1- \Phi(x)$.
Note that $G(r,l) = D_0 + D_1$ where
\be
\label{d0def}
D_0
\ = \
\frac{\sqrt{rl}}{\sqrt{2\pi}}
\int_{-\infty}^{\infty}
e^{-\frac{rlt^2}{2}} \,
\overline{\Phi}^2 \left( \frac{k^2 \tau - rlt}{ \sqrt{k^2 - rl} } \right)
I\{ k^2 \tau - rlt < 1\} \, dt
\ee
and
\be
\label{d1def}
D_1
\ = \
\frac{\sqrt{rl}}{\sqrt{2\pi}}
\int_{-\infty}^{\infty}
e^{-\frac{rlt^2}{2}} \,
\overline{\Phi}^2 \left( \frac{k^2 \tau - rlt}{ \sqrt{k^2 - rl} } \right)
I\{ k^2 \tau - rlt \geq 1\} \, dt .
\ee

\vskip.04in

Consider first the term $D_1$ defined in (\ref{d1def}).
As $rl \neq k^2$ and $k^2\tau-rlt\ge 1$, the normal tail
bound yields
\begin{eqnarray*}
\overline{\Phi} \left( \frac{k^2\tau - rlt}{ \sqrt{k^2 - rl} } \right)
& \leq &
\frac{ \sqrt{k^2 - rl} }{ \sqrt{2 \pi}(k^2 \tau - rlt)} \
\exp\left\{ - \frac{ (k^2 \tau - rlt)^2 }{ 2 (k^2 - rl) } \right\} \nonumber\\
& = &
O(\sqrt{k^2-rl}) \ \exp\left\{-\frac{(k^2\tau-rlt)^2}{2(k^2-rl)}\right\} .
\label{bndphi}
\end{eqnarray*}
Plugging the last expression into (\ref{d1def}), the exponential part of the resulting
integrand is
\[
-\frac{ (k^2 \tau - rlt)^2 }{ (k^2 - rl) } - \frac{rlt^2}{2} ,
\]
which (after lengthy but straightforward algebra) can be expressed as
\[
-k^2\tau+\frac{rl\tau^2}{2}\left(1+\frac{k^2-rl}{k^2+rl}\right)-\frac{rl(k^2+rl)}{2(k^2-rl)}\left((\tau-t)+\tau \left(\frac{k^2-rl}{k^2+rl}\right)\right)^2 .
\]
It then follows that
\begin{eqnarray*}
D_1
& \leq &
O(k^2-rl) \
\exp \left\{ -k^2\tau^2 + \frac{rl\tau^2 }{2}
                  \left( 1 + \frac{k^2-rl}{k^2+rl} \right) \right\}  \nonumber\\
& &
\times \sqrt{\frac{k^2-rl }{ k^2+rl }} \times \int_\infty^\infty \sqrt{\frac{ rl(k^2+rl) }{ k^2-rl }}
\exp
\left\{
-\frac{ rl(k^2+rl) }{2( k^2-rl )}
\left(\tau - t + \frac{ \tau(k^2-rl) }{k^2+rl} \right)^2
\right\}
\it{dt}
\end{eqnarray*}
The term preceding the integral is less than one, and the integral is
equal to one.  Thus $D_1$ is less than the right side of (\ref{grlineq}).

We next consider the term $D_0$ defined in (\ref{d0def}).
Note that $k^2 \tau - rlt < 1$ is equivalent to $t > (k^2 \tau-1) / rl$,
and therefore
\[
D_0
\ \leq \
\int_{ (k^2\tau-1) / rl }^{\infty} \frac{ \sqrt{rl} }{ \sqrt{2\pi} } \, e^{- \frac{rlt^2}{2}} dt
\ = \
\overline{\Phi} \left( \frac{k^2 \tau - 1}{ \sqrt{rl} } \right)
\ \leq \
\frac{k\sqrt{rl}}{\sqrt{2\pi}(k^2\tau-1)}e^{-\frac{(k^2\tau-1)^2}{2rl}-\ln k}\label{l1-2}.
\]
Comparing the last term above with (\ref{grlineq}), it suffices to show that
when $k$ is sufficiently large,
\[
\frac{(k^2 \tau - 1)^2}{2rl} \, + \, \ln k
\ \geq \
\left( k^2-\frac{rl}{2} \right) \, \tau^2
\]
or equivalently
\be
(k^2 - rl)^2 \,  \tau^2 \, - \, 2k^2 \tau \, + \, 1 \, + \, 2rl \ln k \ \geq \ 0.
\label{g3/2q1}
\ee
Suppose first that $rl \geq k^2 - k / \sqrt{\ln k}$. In this case, the left side of
the expression above is at least
\[
-2k^2 \tau \, + \, 1 \, + \, 2rl \ln k
\ \geq \
-2k^2 \tau \, + \, 1 \, + \, 2( k^2 - k / \sqrt{\ln k}) \, \ln k
\  > \ 0
\]
when $k$ is sufficiently large.
Suppose now that $k^2-rl > k / \sqrt{\ln k}$.
As a quadratic function of $\tau$, the left side of
(\ref{g3/2q1}) takes its minimum at $\tau = k^2 / (k^2-rl)^2$, and the
corresponding value is
$rl \, [ -2k^2 + rl + 2(k^2-rl)^2 \ln k ] / (k^2 - rl)^2$.
In this case, the assumption $k^2-rl > k / \sqrt{\ln k}$ implies
\[
-2k^2 \, + \, rl \, + \, 2(k^2 - rl)^2 \, \ln k \ > \ rl \ > \ 0.
\]
This establishes (\ref{g3/2q1}) and complete the proof.

\vskip.2in

\begin{lem}
\label{gauss2}
Let $\tau > 0$ be fixed.  When $k$ is sufficiently large, for every
integer $n$ satisfying the condition
\be
\label{kcondit}
k \ \leq \
\frac{4}{\tau^2} \, \ln n \, - \, \frac{4}{\tau^2} \ln \left(\frac{4}{\tau^2} \ln n\right)
\, - \, \frac{12\ln 2}{\tau^2}
\ee
we have the bound
\[
\frac{ \var \, \Gamma_k(\tau,n)}{(E \, \Gamma_k(\tau,n))^2}
\ \leq \ k^{-2}.
\]
\end{lem}

\vskip.1in

\noindent
{\bf Remark:} For the proof of Theorem \ref{Basic},
it is enough to show that the sum over $k$ of the ratio above is finite,
and for this purpose the upper bound
$k^{-2}$ is sufficient.

\vskip.1in

\noindent{\bf Proof:}
Let $n$ satisfy the condition (\ref{kcondit}).
By Lemma \ref{gauss1}, it suffices to show that
\be
k^4 \,
\sum_{l=1}^k \sum_{r=1}^k \frac{{k \choose l}{n-k \choose k-l}}
       {{n \choose k}} \frac{{k \choose r}{n-k \choose k-r}}{{n \choose k}}
      \, \exp \left\{\frac{rl \tau^2}{2} \left(1+\frac{k^2-rl}{k^2+rl} \right)
     \right\}
\ \leq \
k^{-2}.
\label{g21}
\ee
In order to establish (\ref{g21}), we will show that each term in
the sum is less than $k^{-8}$. To begin, note that
\[
\frac{{k \choose l}{n-k \choose k-l}}{{n \choose k}}
\ \leq \
\frac{ {k \choose l} \, k^l \, (n-k)^{k-l}}{(n-k)^k}
\ = \
{k \choose l} \, k^l(n-k)^{-l},
\]
and that $(n-k)^{-l} = O(n^{-l})$
when $l \le k = O(n^{1/2})$. Thus for some constant $C > 0$,
\[
\frac{{k \choose l}{n-k \choose k-l}}{{n \choose k}}
\frac{{k \choose r}{n-k \choose k-r}}{{n \choose k}}
\ \leq \
C \, {k \choose r}{k \choose l} \, k^{r+l} \, n^{-(r+l)}.
\]
Rewriting (\ref{kcondit}) as $\ln n \geq \frac{\tau^2k}{4} + \ln (\frac{4}{\tau^2}\ln n)+3 \ln 2$
yields the bound
\bea
\lefteqn{
n^{-(r+l)}  \,
   \exp\left\{ \frac{rl\tau^2}{2} \left(1+\frac{k^2-rl}{k^2+rl} \right) \right\} }
   \nonumber \\[.1in]
& \le &
e^{-3 \, (r+l) \, \ln 2} \, \left( \frac{4}{\tau^2} \ln n \right)^{-(r+l)}
  \, \exp \left\{\frac{\tau^2}{2} \left(rl\frac{2k^2}{k^2+rl} \ - \ \frac{k}{2}(r+l)\right)\right\}.\nonumber
\eea
Combining the last three displays, and using the fact that
$k \leq \frac{4}{\tau^2} \ln n$ by assumption,
it suffices to show that
\be
{k \choose r}{k \choose l} \, e^ {-3(r+l) \ln 2} \,
\exp\left\{\frac{\tau^2}{2} \left(rl\frac{2k^2}{k^2+rl}-\frac{k}{2}(r+l)\right)\right\}
\ \leq \
k^{-8}.
\label{targ}
\ee

In order to establish (\ref{targ}), we consider two cases for $r+l$.
Suppose first that $r+l \leq \frac{3k}{4}$.  By elementary
arguments
\[
{k \choose r}{k \choose l} \leq {2k \choose r+l} \leq (2k)^{r+l}
\ \ \mbox{ and } \ \
rl \, \frac{2k^2}{k^2+rl} \ \le \ \frac{(r+l)^2}{4} \, \frac{2k^2}{k^2+rl}
\ \le \ \frac{(r+l)^2}{2} .
\]
It follows from these inequalities that
\bea
\lefteqn{
{k \choose r}{k \choose l}
\, \exp \left\{\frac{\tau^2}{2} \left[ rl \frac{2k^2}{k^2+rl}
\, - \, \frac{k}{2}(r+l) \right] \right\}
} \nonumber \\
& \leq &
\exp \left\{ \frac{\tau^2}{2} \left[ \frac{(r+l)^2}{2} \, - \,
   \frac{k}{2}(r+l) \right] \, + \, (r+l) \, \ln 2k \right\}
   \nonumber \\
& = &
\exp \left\{ \frac{\tau^2(r+l)}{2}
     \left[ \frac{(r+l)}{2} \, - \, \frac{k}{2} \, + \, \frac{2 \ln 2k}{\tau^2} \right]
     \right\} \nonumber \\
& \le &
\exp\left\{ \frac{\tau^2(r+l)}{2}
     \left[ \frac{3k}{8} \, - \, \frac{k}{2} \, + \, \frac{2\ln 2k}{\tau^2} \right] \right\} .
\nonumber
\eea
As the exponent above is negative when $k$ is sufficiently large, (\ref{targ})
follows.
Suppose now that $r+l \geq \frac{3k}{4}$.
From the simple bounds $r+l \geq 2 \sqrt{rl}$ and
$k^2+rl \ge 2\sqrt{k^2rl}$ we find that
\[
rl \, \frac{2k^2}{k^2+rl} \, - \, \frac{k}{2}(r+l)
\ \leq \
\frac{2 rlk^2}{2 \sqrt{k^2 rl}} \, - \, k\sqrt{rl}
\ = \ 0,
\]
and it suffices to bound the initial terms in (\ref{targ}).
But clearly,
\[
{k \choose r}{k \choose l} e^ {-3(r+l)\ln 2}
\ \leq \
2^{2k} \cdot 2^{-\frac{9k}{4}} ,
\]
which is less than $k^{-8}$ when $k$ is sufficiently large.

\vskip.2in

\noindent
{\bf Proof of Theorem \ref{Basic}:}
Proposition \ref{sig-bnd-avg} and the Borel-Cantelli lemma imply that
eventually almost surely
$K_\tau(W_n) \le \lceil s(n,\tau) \rceil+1$.
Thus, we only need to establish an almost sure lower bound on $K_\tau(W_n)$.
To this end, define functions
\[
f(n)
\ = \
\frac{4}{\tau^2} \ln n - \frac{4}{\tau^2} \ln \left(\frac{4}{\tau^2} \ln n\right)
 - \frac{12\ln 2}{\tau^2}
\ \ \mbox{ and } \ \
g(k)
\ = \
\min \{n \geq 1, \lfloor f(n) \rfloor = k\}
\]
for integers $n \geq 1$ and $k \geq 1$, respectively.
It is easy to see that $f(n)$ is strictly increasing for large values
of $n$, and clearly $f(n)$ tends to
infinity as $n$ tends to infinity.
A straightforward argument shows that $g(k)$ has the same properties
Thus for every sufficiently large integer $n$, there
exists a unique integer $k=k(n)$ such that
$g(k) \le n < g(k+1)$.

Fix $m \geq 1$ and consider
the event $A_m$ that for some
$n \geq m$ the random variable $K_\tau(W_n)$ is less than the lower bound
specified in the statement of the theorem.  More precisely, define
\[
A_m \ = \
\bigcup_{n \geq m}
\left\{ K_\tau(W_n) \, \leq \, s(n,\tau) - \frac{12 \ln 2}{\tau^2} - \frac{4}{\tau^2} - 3 \right\}.
\]
To establish the lower bound, it suffices to show that $P(A_m) \to 0$ as $m \to \infty$.
To begin, note that when $m$ is large
\[
A_m
\ \subseteq \
\bigcup_{k \, \geq \, \lfloor f(m)\rfloor} \
\bigcup_{g(k) \, \leq \, n \, < \, g(k+1)}
\left\{ K_{\tau}(W_n) \, \leq \, s(n,\tau) - \frac{12 \ln 2}{\tau^2} - \frac{4}{\tau^2} - 4 \right\} .
\]
Fix $n \geq m$ sufficiently large, and let $k=k(n)$ be the unique
integer such that
$g(k) \le n < g(k+1)$.  The definition of $g(k)$ and the monotonicity
of $f(\cdot)$ ensures that
$k = \lfloor f(g(k)) \rfloor \le f(n) < k+1$.
In conjunction with the definition of $f(n)$ and Lemma \ref{sntau},
this inequality implies that
\bea
1
& = &
k + 1 - k
\ > \
f(n) - \lfloor f(g(k)) \rfloor
\ \ge \
f(n) - f(g(k)) \nonumber \\
& = &
s(n,\tau) - s(g(k),\tau) + o(1), \nonumber
\eea
and therefore $s(n,\tau)< s(g(k),\tau)+1+o(1)$.  Define
\[
r(k) \ = \ \left\lfloor s(g(k),\tau) - \frac{12 \ln 2}{\tau^2} - \frac{4}{\tau^2} \right\rfloor .
\]
From the bound on $s(n,\tau)$ above
and the fact that $K_\tau (W_{g(k)}) \leq K_{\tau}(W_n)$, we have
\begin{eqnarray*}
\left\{ K_{\tau}(W_n) \, \leq \, s(n,\tau) - \frac{12 \ln 2}{\tau^2} - \frac{4}{\tau^2} - 3 \right\}
& \subseteq &
\left\{ K_{\tau}(W_{g(k)}) \, \leq \, r(k) - 1 + o(1) \right\} \\
& \subseteq &
\left\{ K_{\tau}(W_{g(k)}) \, \leq \, r(k) - 1 \right\} ,
\end{eqnarray*}
where the last relation makes use of the fact that $K_\tau$ and
$r(k)$ are integers.  Thus we find that
\[
A_m
\ \subseteq \
\bigcup_{k \, \ge \, \lfloor f(m) \rfloor}
\left\{ K_{\tau}(W_{g(k)})
         \ \leq \ r(k) - 1 \right\} .
\]
Consider the events above.  For fixed $k$,
\be
\label{gk}
P( K_{\tau}(W_{g(k)}) \leq r(k) -1 )
\ = \
P( \Gamma_{r(k)}(\tau, g(k)) = 0 )
\ \leq \
\frac{ \mbox{Var} \, \Gamma_{r(k)} (\tau, g(k))}{ (E \, \Gamma_{r(k)}(\tau,g(k)))^2}
\ee
where we have used the fact that for a non-negative integer-valued random
variable $X$
\[
P( X = 0 ) \ \leq \ P( | X - EX | \geq EX ) \ \leq \
\frac{\mbox{Var} X}{ (E X)^2}
\]
by Chebyshev's inequality.  As $r(k) \leq f(g(k))$, Lemma \ref{gauss2}
ensures that the final term in (\ref{gk}) is less than $k^{-2}$, and
the Borel-Cantelli lemma then implies that $P(A_m) \rightarrow 0$ as
$m \rightarrow \infty$.  This completes the proof of Theorem \ref{Basic}.

\vskip.3in

\noindent
{\bf\Large  Acknowledgements} \\
The authors would like to thank Andrey Shabalin for his
assistance with the simulation results in Section{}, and
for his help in clarifying the connections between the work
described here and results in random matrix theory.
We would also like to thank John Hartigan for pointing out the use
of the Gaussian comparison principle as an alternative way
of obtaining the bounds of Proposition \ref{sig-bnd-avg}.
The work presented in this paper was supported in part
by NSF grants DMS 0406361 and DMS 0907177.


\begin{thebibliography}{9}
\csname bibmessage\endcsname

\bibitem{alon04}
\textsc{Alon, N.} and
\textsc{Naor, A.} (2006).
Approximating the Cut-Norm via Grothendieck's Inequality.
\textit{SIAM Journal of Computating}
\textbf{35:4}, 787--803.

\bibitem{gold90}
\textsc{Arratia, R.}, \textsc{Goldstein, L.} and
\textsc{Gordon, L.} (1990).
Poisson Approximation and the Chen-Stein Method.
\textit{Statistical Science}
\textbf{5:4}, 403--424.

\bibitem{bol2}
\textsc{Bollob\'as, B.} and
\textsc{Erd\H{o}s, P.} (1976).
Cliques in Random Graphs.
\textit{Mathematical Proceedings of the Cambridge Philosophical Society}
\textbf{80}, 419--427.

\bibitem{cheng00}
\textsc{Cheng, Y.} and
\textsc{Church, G. M.} (2000).
Biclustering of Expression Data.
\textit{Proceedings of the 8th International Conference on Intelligent Systems for Molecular Biology} 93--103.

\bibitem{dhi1}
\textsc{Dawande, M.}, \textsc{Keskinocak, P.}, \textsc{Swaminathan, J. M.} and
\textsc{Tayur, S.} (2001).
On Bipartite and Multipartite Clique Problems.
\textit{Journal of Algorithms}
\textbf{41:2}, 388--403.

\bibitem{garber2001dge}
\textsc{Garber, M.E.}, \textsc{Troyanskaya, O.G.}, \textsc{Schluens, K.}, \textsc{Petersen, S.}, \textsc{Thaesler, Z.},
\textsc{Pacyna-Gengelbach, M.}, \textsc{van de Rijn, M.}, \textsc{Rosen, G.D.}, \textsc{Perou, C.M.}, \textsc{Whyte, R.I.}, \textsc{Altman, R.B.}, \textsc{Brown, P.O.},
\textsc{Botstein, D.} and
\textsc{Petersen, I.} (2001).
Diversity of Gene Expression in Adenocarcinoma of the Lung.
\textit{Proceedings of the National Academy of Sciences}
\textbf{98:24}, 13784--13789.

\bibitem{geman}
\textsc{Geman, S.} (1980).
A Limit Theorem for the Norm of Random Matrices.
\textit{Annals of Probability}
\textbf{8:2}, 252--261.

\bibitem{golub1999mcc}
\textsc{Golub, T.R.}, \textsc{Slonim, D.K.}, \textsc{Tamayo, P.}, \textsc{Huard, C.}, \textsc{Gaasenbeek, M.},
\textsc{Mesirov, J.P.}, \textsc{Coller, H.}, \textsc{Loh, M.L.}, \textsc{Downing, J.R.}, \textsc{Caligiuri, M.A.}, \textsc{Bloomfield, C.D.} and
\textsc{Lander, E.S.} (1999).
Molecular Classification of Cancer: Class Discovery and Class Prediction by Gene Expression Monitoring.
\textit{Science}
\textbf{286:5439}, 531--537.


\bibitem{matrixhj}
\textsc{Horn, R.A.} and
\textsc{Johnson, C.R.} (1985).
Matrix Analysis.
\textit{Cambridge University Press}.

\bibitem{jiang2004cag}
\textsc{Jiang, D.}, \textsc{Tang, C.} and
\textsc{Zhang, A.} (2004).
Cluster Analysis for Gene Expression Data: A Survey.
\textit{IEEE Transactions on Knowledge and Data Engineering}
\textbf{16:11}, 1370--1386.

\bibitem{Gar2}
\textsc{Koyuturk, M.}, \textsc{Szpankowski, W.} and
\textsc{Grama, A.} (2004).
Biclustering Gene-Feature Matrices for Statistically Significant Dense Patterns.
\textit{Proceedings of the 2004 IEEE Computational Systems Bioinformatics Conference} 480--484.

\bibitem{owen02}
\textsc{Lazzeroni, L.} and \textsc{Owen, A.} (2000).
Plaid Models for Gene Expression Data.
\textit{Statistica Sinica}
\textbf{12}, 61--86.

\bibitem{madeira2004bab}
\textsc{Madeira, S.C.} and
\textsc{Oliveira, A.L.} (2004).
Biclustering Algorithms for Biological Data Analysis: A Survey.
\textit{IEEE/ACM Transactions on Computational Biology and Bioinformatics}
\textbf{1:1}, 24--45.

\bibitem{Mat1}
\textsc{Matula, D.} (1976).
The Largest Clique Size in A Random Graph.
\textit{Southern Methodist University Technical Report},
CS 7608.


\bibitem{parsons2004sch}
\textsc{Parsons, L.}, \textsc{Haque, E.} and
\textsc{Liu, H.} (2004).
Subspace Clustering for High Dimensional Data: A Review.
\textit{ACM SIGKDD Explorations Newsletter}
\textbf{6:1}, 90--105.

\bibitem{paul}
\textsc{Paul, D.} (2007).
Asymptotics of Sample Eigenstructure for A Large Dimensional Spike Covariance Model.
\textit{Statistica Sinica}
\textbf{17}, 1617--1642.

\bibitem{perou2000mph}
\textsc{Perou, C.M.}, \textsc{S{\o}rlie, T. }, \textsc{Eisen, M.B.}, \textsc{van de Rijn, M. }, \textsc{Jeffrey, S.S.}, \textsc{Rees, C.A.}, \textsc{Pollack, J.R.},
\textsc{Ross, D.T.}, \textsc{Johnsen, H.}, \textsc{Akslen, L.A.}, \textsc{Fluge, {\O}.}, \textsc{Pergamenschikov, A.}, \textsc{Williams, C.}, \textsc{Zhu, S.X.},
\textsc{L{\o}nning, P.E.}, \textsc{B{\o}rresen-Dale, A.L.}, \textsc{Brown, P.O.} and \textsc{Botstein, D.}(2000).
Molecular Portraits of Human Breast Tumors.
\textit{Nature}
\textbf{406}, 747--752.

\bibitem{ShabalinandNobel}
\textsc{Shabalin, A.A.}, \textsc{Weigman, V.J.}, \textsc{Perou, C.M.} and \textsc{Nobel, A.B.}
(2009).
Finding Large Average Submatrices in High Dimensional Data.
\textit{Annals of Applied Statistics}
\textbf{3}, 985--1012.


\bibitem{slep62}
\textsc{Slepian, D.} (1962).
The One-sided Barrier Problem for Gaussian Noise.
\textit{Bell System Technical Journal}
\textbf{41}, 463--501.

\bibitem{sorlie2001gep}
\textsc{S{\o}rlie, T.}, \textsc{Perou, C.M.}, \textsc{Tibshirani, R.}, \textsc{Aas, T.}, \textsc{Geisler, S.}, \textsc{Johnsen, H.}, \textsc{Hastie, T.},
\textsc{Eisen, M.B.}, \textsc{van de Rijn, M.}, \textsc{Jeffrey, S.S.}, \textsc{Thorsen, T.}, \textsc{Quist, H.}, \textsc{Matese, J.C.}, \textsc{Brown, P.O.},
\textsc{Botstein, D.}, \textsc{L{\o}nning, P.E.} and \textsc{B{\o}rresen-Dale, A.} (2001).
Gene Expression Patterns of Breast Carcinomas Distinguish Tumor Subclasses with Clinical Implications.
\textit{Proceedings of the National Academy of Sciences}
\textbf{98:19}, 10869--10874.

\bibitem{sorlie2003rob}
\textsc{S{\o}rlie, T.}, \textsc{Tibshirani, R.}, \textsc{Parker, J.}, \textsc{Hastie, T.}, \textsc{Marron, J.S.}, \textsc{Nobel, A.B.},
\textsc{Deng, S.}, \textsc{Johnsen, H.}, \textsc{Pesich, R.}, \textsc{Geisler, S.}, \textsc{Demeter, J.}, \textsc{Perou, C. M.},
\textsc{L{\o}nning, P.E.}, \textsc{Brown, P.O.}, \textsc{B{\o}rresen-Dale, A.} and \textsc{Botstein, D.}
(2003).
Repeated Observation of Breast Tumor Subtypes in Independent Gene Expression Data Sets.
\textit{Proceedings of the National Academy of Sciences}
\textbf{100:14}, 8418--8423.

\bibitem{sn06}
\textsc{Sun, X.} and
\textsc{Nobel, A.B.} (2006).
Significance and Recovery of Block Structures in Binary Matrices with Noise.
\textit{Proceedings of the 19th Conference on Computational Learning Theory}, 109--122.

\bibitem{sn08}
\textsc{Sun, X.} and
\textsc{Nobel, A.B.} (2008).
On the Size and Recovery of Submatrices of Ones in a Random Binary Matrix.
\textit{Journal of Machine Learning Research}
\textbf{9}, 2431--2453.

\bibitem{tany1}
\textsc{Tanay, A.}, \textsc{Sharan, R.} and
\textsc{Shamir, R.} (2002).
Discovering Statistically Significant Biclusters in Gene Expression Data.
\textit{Bioinformatics}
\textbf{18:1}, 136--144.

\bibitem{weigelt2005mpa}
\textsc{Weigelt, B.}, \textsc{Hu, Z.}, \textsc{He, X.}, \textsc{Livasy, C.}, \textsc{Carey, L.A.}, \textsc{Ewend, M.G.}, \textsc{Glas, A.M.}, \textsc{Perou, C.M.}
and \textsc{van't Veer, L.~J.} (2005).
Molecular Portraits and 70-Gene Prognosis Signature are Preserved throughout the Metastatic Process of Breast Cancer.
\textit{Cancer Research}
\textbf{65:20}, 9155--9158.
\end{thebibliography}
\end{document}